\input amstex
\documentstyle{amsppt}
\magnification=1200
\NoBlackBoxes
\hsize=5.6in
\vsize=7.3in
\topmatter
\title Distributional Limits for the Symmetric Exclusion Process\endtitle
\rightheadtext {Symmetric Exclusion Processes}
\date October 18, 2007\enddate
\author  Thomas M. Liggett\endauthor
\thanks Research supported in part by
NSF Grant DMS-0301795.\endthanks
\affil University of California, Los Angeles\endaffil
\abstract Strong negative dependence properties have recently been proved for the
symmetric exclusion process. In this paper, we apply these results to prove convergence
to the Poisson and normal distributions for various functionals of the process.
\endabstract 
\endtopmatter

\heading 1. Introduction\endheading
The symmetric exclusion process on the countable set $S$ is the Markov process $\eta_t$ on $\{0,1\}^S$
with formal generator
$$\Cal Lf(\eta)=\sum_{\eta(x)=1,\eta(y)=0}p(x,y)[f(\eta_{x,y})-f(\eta)],$$
where $\eta_{x,y}$ is the configuration obtained from $\eta$ by interchanging the
coordinates $\eta(x)$ and $\eta(y)$. Here $p(x,y)=p(y,x)$ are the transition probabilities
for a symmetric, irreducible, Markov chain on $S$. 
For background on this process, see Chapter VIII of Liggett (1985). 

Many limit theorems of various types have
been proved for this process. Examples are the central limit theorems for a tagged particle
and for the flux in one dimensional systems by
Arratia (1983), Kipnis and Varadhan (1986), De Masi and Ferrari (2002), Jara and Landim (2006), and Peligrad
and Sethuraman (2008). In this paper, we focus on limit theorems that can now be proved using the recently
obtained strong negative dependence properties of the symmetric exclusion process.

A probability measure $\mu$ on $\{0,1\}^S$ is said to be negatively associated if 
$$\int fgd\mu\leq\int fd\mu\int gd\mu$$
for all increasing continuous functions $f,g$ on $\{0,1\}^S$ that depend on disjoint sets of
coordinates. Theorem 5.2 of
Borcea, Br\"and\'en and  Liggett (2008) asserts that if the initial distribution of the symmetric exclusion
process
$\eta_t$ is a product measure, then the distribution of $\eta_t$ is negatively associated for all $t>0$. 
In fact, by Proposition 5.1 of that paper, it has a stronger and even more useful property, known as strong
Rayleigh.

Limit theorems for negatively associated random variables have been proved by a number of
authors -- see Barbour, Holst and Janson (1992), Newman (1984) and Roussas (1994), for example. In the case of
convergence to the normal law, none of these results
quite fit our setting. In our situation, there is generally no translation invariance in the
covariance structure, and the sum of off-diagonal covariances is often not ``little o" of the
sum of variances. However, we will see in
Section 2 that the strong Rayleigh property makes it quite easy to prove convergence to the Poisson and
Gaussian laws, given estimates of variances and covariances. Therefore, we will not need to
use these earlier results.

The first situation we will consider involves the extremal stationary distributions for the process.
We recall their description -- see Chapter VIII of Liggett (1985). Let
$$\Cal H=\bigg\{\alpha: S\rightarrow [0,1], \sum_yp(x,y)\alpha(y)=\alpha(x)\ \forall x\bigg\},$$
and for $\alpha\in\Cal H$, let $\nu_{\alpha}$ be the product measure with marginals
$\nu\{\eta:\eta(x)=1\}=\alpha(x).$ Then the limiting distribution as $t\rightarrow\infty$ of the
process at time $t$ exists if the initial distribution is $\nu_{\alpha}$; call it $\mu_{\alpha}$. 
The result is that the extremal stationary distributions are exactly $\{\mu_{\alpha},\alpha\in \Cal H\}$.
If $\alpha$ is constant, then $\mu_{\alpha}=\nu_{\alpha}$ so we are really interested in nonconstant
$\alpha$'s, in which case very little is known about the corresponding stationary
distributions other than the marginals -- $\mu_{\alpha}\{\eta:\eta(x)=1\}
=\alpha(x)$. If $p(x,y)$ are the
transition probabilities for simple random walk on a homogeneous tree, for example, there are many such
nonconstant $\alpha$'s, and therefore there are many extremal stationary distributions that are not of
product form. We now know from the results in Borcea, Br\"and\'en and 
Liggett (2008) that $\nu_{\alpha}$ is negatively associated -- and in fact strong Rayleigh -- for each
$\alpha\in
\Cal H$.

We will use the following notation. For $n\geq 1$, $p^{(n)}(x,y)$ are the $n$-step transition
probabilities for the Markov chain with transition probabilities $p(x,y)$, and for $t>0$,
$$p_t(x,y)=e^{-t}\sum_{n=0}^{\infty}\frac{t^n}{n!}p^{(n)}(x,y)$$
are the transition probabilities for the corresponding continuous time chain $X_t$. The Green
function is given by
$$G(x,y)=\sum_{n=0}^{\infty}p^{(n)}(x,y)=\int_0^{\infty}p_t(x,y)dt.$$
The Dirichlet sum of an $\alpha\in\Cal H$ is defined by 
$$\Phi(\alpha)=\sum_{x, y}p(x,y)[\alpha(y)-\alpha(x)]^2.$$
This quantity is finite for many, but not all, elements of $\Cal H$ if $S$ is a regular tree,
for example. A construction of a class of infinite graphs with only one end that support
nonconstant $\alpha\in\Cal H$ with $\Phi(\alpha)<\infty$ is constructed in Cartwright and Woess
(1992). In this context, we have the following results. We use $\Rightarrow$ to denote convergence
in distribution.

\proclaim{Theorem 1}Suppose $\alpha\in\Cal H$ and $\Phi(\alpha)<\infty.$ If $S_n\subset S$ 
satisfy

(a) $\lim_{n\rightarrow\infty}\sup_{x\in S_n}\alpha(x)=0,$ $\lim_{n\rightarrow\infty}\sum_{x\in
S_n}\alpha(x)=\lambda<\infty,$ 

\noindent and

(b) $\sup_{x,n}\sum_{y\in S_n}G(x,y)<\infty,$

\noindent then under $\mu_{\alpha},$
$$\sum_{x\in S_n}\eta(x)\Rightarrow Poisson(\lambda).$$\endproclaim

\proclaim{Theorem 2} Suppose $\alpha\in\Cal H$ and $\Phi(\alpha)<\infty.$ If $S_n\subset S$  satisfy
$$\lim_{n\rightarrow\infty}\sum_{x\in S_n}\alpha(x)[1-\alpha(x)]=\infty$$

\noindent and

$$\liminf_{n\rightarrow\infty}\frac{\sum_{x\in S_n}\alpha(x)[1-\alpha(x)]}{\sup_x\sum_{y\in S_n}G(x,y)}
>\Phi(\alpha)\tag 1.1$$

\noindent then under $\mu_{\alpha}$,
$$\frac{\sum_{x\in S_n}\eta(x)-\sum_{x\in S_n}\alpha(x)}{[Var_{\mu_{\alpha}}(\sum_{x\in
S_n}\eta(x))]^{1/2}}\Rightarrow  N(0,1).$$
Furthermore,
$$Var_{\mu_{\alpha}}\bigg(\sum_{x\in S_n}\eta(x)\bigg)\bigg/\sum_{x\in S_n}\alpha(x)[1-\alpha(x)],\tag 1.2$$
(which is at most one) is bounded below by a positive constant. If the left side of (1.1) is
infinite, then the limit of (1.2) as $n\rightarrow\infty$ is 1.
\endproclaim

Theorems 1 and 2 will be proved in Section 3, after deriving limit theorems for general strong
Rayleigh Bernoulli random variables in Section 2. 

As an example of the application of Theorems 1 and 2, 
let $S$ be the binary tree, and let the chain have nearest neighbor jumps with probability 1/3 each.
Write $S=L\cup R$ where $L,R$ are defined as follows: An  basis edge is fixed, and its endpoints are called
left and right respectively. Then $L$ is the set of vertices that are closer to the left vertex than
to the right, and $R=S\backslash L$. Each vertex $x$ in $S$ is assigned a level $l(x)\geq 0$, which is
the distance from $x$ to the closer of the two endpoints of the basis edge. Then $G(x,y)=2^{-d(x,y)+1}$,
where $d(x,y)$ is the distance between $x$ and $y$. Therefore, $\sup_x\sum_{y:l(y)<n}G(x,y)=3n$,
$\sup_x\sum_{y\in L:l(y)<n}G(x,y)=2n$ and $\sup_x\sum_{y\in L:l(y)=n}G(x,y)=3-2^{-n}$.

For $0\leq \lambda,\rho\leq 1$,
let $\alpha\in \Cal H$ be defined by
$$\alpha(x)=\cases\lambda+\frac{\rho-\lambda}{3\cdot 2^{l(x)}}&\quad\text{if }x\in L,\\
\rho+\frac{\lambda-\rho}{3\cdot 2^{l(x)}}&\quad\text{if }x\in R,\endcases$$
and put $\mu=\mu_{\alpha}$.
Then $\Phi(\alpha)=2(\rho-\lambda)^2/9$,
$$\gather E^{\mu}\sum_{x:l(x)<n}\eta(x)=(\lambda+\rho)(2^n-1),\quad\text{and}\\\sum_{x:l(x)<
n}Var_{\mu}\eta(x)= (2^n-1)[\rho(1-\rho)+\lambda(1-\lambda)]+(\lambda-\rho)^2\bigg[\frac 23 n-\frac
49(1-2^{-n})\bigg].\endgather$$
It follows that for $S_n=\{x:l(x)<n\}$, the left side of (1.1) is infinite if 
$$\rho(1-\rho)+\lambda(1-\lambda)>0.$$
 Therefore, for all choices of $\lambda,\rho$, Theorem 2 implies that
$$\frac{\sum_{x:l(x)<n}\eta(x)-(\lambda+\rho)2^n}{\sqrt{2^n}}\Rightarrow
N(0,\rho(1-\rho)+\lambda(1-\lambda)).\tag 1.3$$

Next take $\lambda=0, \rho=1$, in which case (1.3) has little content. Then
Theorem 1 gives
$$\sum_{x\in L:l(x)=n}\eta(x)\Rightarrow Poisson (1/3).$$
If $S_n=\{x\in L:l(x)<n\}$, the left and right sides of (1.1) are $\frac 16$ and $\frac 29$ respectively, so
(1.1) does not hold in this case. Nevertheless, we will see at the end of Section 3 that
more careful estimates imply that 
$$\frac{\sum_{x\in L:l(x)<n}\eta(x)-\frac n3}{\sigma_n}\Rightarrow  N(0,1),\tag 1.4$$
where $\sigma_n^2/n$ is asymptotically between $\frac {23}{189}$ and $\frac 13$.

The next situation we consider was proposed by Pemantle (2000) as an application of the then hoped for
negative dependence properties of the symmetric exclusion process. Now $S=Z^1$, and $p(x,y)=
p(y-x)$, with $\sum_x|x|p(x)<\infty$. For the initial configuration, we take
$$\eta(x)=\cases 1\quad&\text{if }x\leq 0,\\0\quad&\text{if }x> 0.\endcases$$
Since $\sum_{x<0<y}p(x,y)<\infty$, $W_t=\sum_{x>0}\eta_t(x)<\infty$ a.s. for all $t$.

\proclaim{Theorem 3} Suppose $\sigma^2=\sum_nn^2p(n)<\infty$. Then
$$\frac {W_t-EW_t}{[Var(W_t)]^{1/2}}\Rightarrow N(0,1)\tag 1.5$$
as $t\rightarrow\infty$. Furthermore, 
$$\lim_{t\rightarrow\infty}\frac {EW_t}{\sqrt t}=\frac{\sigma}{\sqrt{2\pi}},$$
and
$$\frac{Var(W_t)}{t^{1/2}} $$
is bounded above and below by positive constants.
\endproclaim
Theorem 3 will be proved in Section 4. It seems likely that if the distribution $p(\cdot)$
is in the domain of normal attraction of a (symmetric) stable law of index $\alpha\in (1,2)$,
then (1.5) holds with $Var(W_t)$ of order $t^{1/\alpha}$, but this has not been checked.

The major result from Borcea, Bra\"nden and Liggett (2008) that we use in this paper is Proposition
5.1, which asserts that the strong Rayleigh property is preserved by the evolution of a symmetric
exclusion process. The proof given there uses results from earlier papers, including Obreschkoff's
Theorem, which are not part of the toolkit of a typical probabilist. In Section 5, we present an
elementary proof of that result in order to make the present paper essentially self-contained.

\heading 2. Limit theorems for strong Rayleigh measures\endheading

Consider a probability measure $\mu$ on $\{0,1\}^n$. Its generating
polynomial is defined for $z\in \Cal C^n$ by
$$Q(z)=Q(z_1,...,z_n)=E^{\mu}\prod_{i=1}^n z_i^{\eta(i)}.$$
The measure $\mu$ is said to be Rayleigh if
$$\frac{\partial Q}{\partial z_i}(z)\frac{\partial Q}{\partial z_j}(z)\geq Q(z)\frac{\partial^2 Q}
{\partial z_i\partial z_j}(z),\quad i\neq j\tag 2.1$$
for all $z\in R_+^n$, and to be strong Rayleigh if (2.1) holds for all $z\in R^n$. Note that
when $z=(1,1,...,1)$, (2.1) says that $\eta(i)$ and $\eta(j)$ are negatively correlated under $\mu$.
By Theorem 4.9 of Br\"anden, Borcea and Liggett (2008), the strong Rayleigh property implies
negative association.

Br\"and\'en (2007) proved that the strong Rayleigh property is equivalent to the following
property, which is known as stability: $Q(z)\neq 0$ if $z_i$ has strictly positive imaginary
part for each $i$. This is the key to the following representation.

\proclaim{Proposition 2.1} Suppose $\mu$ is strong Rayleigh. Then there exist independent Bernoulli
random variables $\zeta_i$ with parameters $p_i$ so that the distribution of $\sum_i\eta(i)$ under $\mu$
is the same as that of $\sum_i\zeta_i$.\endproclaim

\demo{Proof} Setting $z_i\equiv w$ in the expression for $Q$, we see that the polynomial
in one variable
$$Q^*(w)=Q(w,w,...,w)=E^{\mu}w^{\sum_i\eta(i)}$$
has no roots with positive imaginary part, and therefore all of its roots are real. Since
$Q^*(1)=1$, it follows that $Q^*$ can be written in the form
$$Q^*(w)=\prod_i(w+w_i)\bigg/\prod_i(1+w_i)$$
where the $w_i$ are real. Since $Q^*(-w_i)=0$, we see that $w_i\geq 0$ for each $i$. 
Letting $p_i=1/(1+w_i)$, this becomes
$$Q^*(w)=\prod_i\big[p_iw+(1-p_i)\big],$$
which is the generating polynomial for $\sum_i\zeta_i$. \enddemo

Using this result, it is easy to extend the classical limit theorems to triangular
arrays of strong Rayleigh random variables.

\proclaim{Proposition 2.2} Suppose the Bernoulli random variables $\{\eta_n(x)\}$ are
strong Rayleigh for each $n$.

(a) If $\lim_{n\rightarrow\infty}\sum_xE\eta_n(x)=\lambda,\ \lim_{n\rightarrow\infty}\sum_x[E\eta_n(x)]^2=0,$
and 
$$\lim_{n\rightarrow\infty}\sum_{x\neq y}Cov(\eta_n(x),\eta_n(y))=0,$$
then 
$$\sum_x\eta_n(x)\Rightarrow Poisson (\lambda).$$

(b) If $\lim_{n\rightarrow\infty}Var(\sum_x\eta_n(x))=\infty$, then

$$\frac{\sum_x\eta_n(x)-E\sum_x\eta_n(x)}{\sqrt{Var(\sum_x\eta_n(x))}}\Rightarrow N(0,1).$$

\endproclaim

\demo{Proof} Using Proposition 2.1, let $\zeta_{n,i}$ be Bernoulli random variables that are
independent in $i$ for each $n$, and have the property that
$$\sum_x\eta_n(x)\quad\text{ and }\quad\sum_i\zeta_{n,i}$$
have the same distribution for each $n$. It suffices to show that the conditions for Poisson
and normal convergence hold for the array $\zeta_{n,i}$. But this is immediate from the
assumptions and the following identities:
$$\sum_xE\eta_n(x)=\sum_iE\zeta_{n,i}$$
and
$$\align Var\bigg(\sum_x\eta_n(x)\bigg)=&\sum_xVar(\eta_n(x))+\sum_{x\neq y}Cov(\eta_n(x),\eta_n(y))\\=&
Var\bigg(\sum_i\zeta_{n,i}\bigg)=\sum_iVar(\zeta_{n,i}).\endalign$$
\enddemo

\noindent {\bf Remark.} Proposition 2.2(a) is a consequence of Theorem 11 in Newman (1984) or
Theorem 3D in Barbour, Holst and Janson (1992). The latter also gives bounds on the total
variation distance from the Poisson distribution. We note that
under the strong Rayleigh assumption we are making, the proof is very simple.

\heading 3. The stationary distributions\endheading

We begin the proofs of Theorems 1 and 2 by obtaining a bound on the covariances for the measure
$\mu=\mu_{\alpha}$. Let $U$ and $U(t)$ be the generator and semigroup for the motion of two independent copies
of  the Markov chain with transition probabilities $p_t(x,y)$, and $V$ and $V(t)$ be the generator
and semigroup for the motion of two copies of that Markov chain with the exclusion interaction. Then
by (1.28) and (1.29) on page 373 of Liggett (1985),
$$E^{\mu}\eta(x)\eta(y)=\lim_{t\rightarrow\infty}V(t)f(x,y),\quad x\neq y,$$
where $f(x,y)=\alpha(x)\alpha(y)$. Since $U(t)f(x,y)=f(x,y)$, the integration by parts formula
gives
$$\aligned -Cov_{\mu}(\eta(x),\eta(y))=&f(x,y)-\lim_{t\rightarrow\infty} V(t)f(x,y)
\\=&\lim_{t\rightarrow\infty}\bigg[U(t)f(x,y)-V(t)f(x,y)\bigg]\\=&\lim_{t\rightarrow\infty}\int_0^{t}
V(s)[U-V]U(t-s)f(x,y)ds\\=&\int_0^{\infty}
V(s)[U-V]f(x,y)ds\endaligned\tag 3.1$$
for $x\neq y$. From page 366 of Liggett (1985), we see that for $x\neq y$,
$$(U-V)f(x,y)=p(x,y)[f(x,x)+f(y,y)-2f(x,y)]=p(x,y)[\alpha(x)-\alpha(y)]^2.\tag 3.2$$
Define $\Delta(x,y)=p(x,y)[\alpha(x)-\alpha(y)]^2$ for all $x,y$.

Let $g_n(x,y)=1_{S_n\times S_n}(x,y).$ This is a positive definite function. Therefore, using the
symmetry of $V(s)$, the fact that $\Delta(x,y)=0$ if $x=y$,  and Proposition 1.7 on page 366 of Liggett (1985),
we see that
$$\aligned -\sum_{x,y\in S_n;x\neq y}Cov_{\mu}(\eta(x),\eta(y))=&\int_0^{\infty}\sum_{x\neq
y}g_n(x,y)V(s)\Delta(x,y)ds\\=&\int_0^{\infty}\sum_{x,y}\Delta(x,y)V(s)g_n(x,y)ds\\
\leq& \int_0^{\infty}\sum_{x,y}\Delta(x,y)U(s)g_n(x,y)ds\\=&\sum_{x,y}\Delta(x,y)\int_0^{\infty}
P^x(X_s\in S_n)P^y(X_s\in S_n)ds,\endaligned\tag 3.3$$
where $X_s$ is the Markov chain with transition probabilities $p_s(x,y)$.

\demo{Proof of Theorem 1} Given the strong Rayleigh property of $\mu$ and Proposition 2.2(a), it suffices to
check that
$$\lim_{n\rightarrow\infty}\sum_{x,y\in S_n;x\neq y}Cov_{\mu}(\eta(x),\eta(y))=0,$$
and therefore we need to check that the right side of (3.3) tends to zero. To do so, note that
$$\int_0^{\infty}P^x(X_s\in S_n)ds=\sum_{u\in S_n}G(x,u)\leq C_1\tag 3.4$$
for some constant $C_1$ by assumption (b). Next, since
$$\sup_{s>0}e^{-s}\frac{s^k}{k!}=e^{-k}\frac{k^k}{k!}\leq \frac {C_2}{\sqrt k}$$
for some constant $C_2$, we have for any $N\geq 1$,
$$P^y(X_s\in S_n)\leq \sum_{k=0}^N\sum_{v\in S_n}p^{(k)}(y,v)+\frac {C_1C_2}
{\sqrt N},$$
so that (using the first part of assumption (a)) for each $y$
$$\lim_{n\rightarrow\infty}\sup_{s>0}P^y(X_s\in S_n)=0.\tag 3.5$$
It follows that the right side of (3.3) tends to zero by (3.4), (3.5), the finiteness of $\Phi(\alpha)$ and the
dominated convergence theorem.\enddemo

\demo{Proof of Theorem 2} By Proposition 2.2(b), we need to show that
$$\lim_{n\rightarrow\infty}Var_{\mu}\bigg(\sum_{x\in S_n}\eta(x)\bigg)=\infty.\tag 3.6$$
By (3.3),
$$-\sum_{x,y\in S_n;x\neq y}Cov_{\mu}(\eta(x),\eta(y))\leq \Phi(\alpha)\sup_x\sum_{u\in S_n}G(x,u).$$
Therefore, (3.6) follows from the assumptions in the theorem.
\enddemo

We now return to the example discussed in the introduction -- the simple random walk on the binary
tree with
$$\alpha(x)=\cases \frac 1{3\cdot 2^{l(x)}}\quad&\text{if }x\in L,\\
1-\frac 1{3\cdot 2^{l(x)}}\quad&\text{if }x\in R,\endcases$$
and $S_n=\{x\in L:l(x)<n\}$. We will see that Proposition 2.2(b) applies, even though (1.1) fails.
In order to do so, we will use the structure of the problem to estimate the right side of (3.3)
more carefully. 

First note that $l(X_t)$ is a Markov chain on the nonnegative integers with drift $\frac 13$.
Therefore 
$$\frac{l(X_t)}t\rightarrow\frac 13\quad a.s.\tag 3.7$$
by the strong law of large numbers. Secondly, if $\beta,\gamma>0$ are chosen so that
$$e^{-\gamma}+e^{+\gamma}=3(1-\beta),$$
then
$$e^{\beta t-\gamma l(X_t)}$$
is a supermartingale, so that
$$1\geq E^xe^{\beta nt-\gamma l(X_{nt})}\geq e^{\beta n t-\gamma n}P^x(l(X_{nt})<n).$$
This implies
$$P^x(l(X_{nt})<n)\leq e^{n(\gamma-\beta t)}\leq e^{\gamma-\beta t}$$
if $n\geq 1$ and $t\geq \gamma/\beta$, which provides the domination in the following
computation.
$$\aligned\lim_{n\rightarrow\infty}\frac 1n\sum_{x,y}\Delta(x,y)&\int_0^{\infty}
P^x(X_s\in S_n)P^y(X_s\in S_n)ds\\=&\lim_{n\rightarrow\infty}\sum_{x,y}\Delta(x,y)\int_0^{\infty}
P^x(X_{nt}\in S_n)P^y(X_{nt}\in S_n)dt\\=&3\sum_{x,y}\Delta(x,y)[1-\alpha(x)][1-\alpha(y)].
\endaligned\tag 3.8$$
For the final step above, note that the integrand in the middle line tends to zero if $t>3$
by (3.7), while if $t<3$,
$$\lim_{n\rightarrow\infty}P^x(X_{nt}\in S_n)=P^x(X_s\in L\text{ eventually})=1-\alpha(x).$$

Next, we compute the right side of (3.8).
$$\aligned &3\sum_{x,y}\Delta(x,y)[1-\alpha(x)][1-\alpha(y)]=\sum_{d(x,y)=1}[\alpha(x)-\alpha(y)]^2
[1-\alpha(x)][1-\alpha(y)]\\=&2\bigg[\frac 2{3^4}+\sum_{n=1}^{\infty}2^n\frac 1{(3\cdot 2^n)^2}
\frac 1{3\cdot 2^n}\frac 1{3\cdot 2^{n-1}}+\sum_{n=1}^{\infty}2^n\frac 1{(3\cdot 2^n)^2}
\bigg(1-\frac 1{3\cdot 2^n}\bigg)\bigg(1-\frac 1{3\cdot 2^{n-1}}\bigg)\bigg]\\=&\frac{40}{189}<\frac 13.
\endaligned$$
Combining this with (3.3) and (3.8), we see that
$$\liminf_{n\rightarrow\infty}\frac{Var\bigg(\sum_{x\in S_n}\eta(x)\bigg)}n>0,$$
so that Proposition 2.2(b) gives (1.4).

\heading 4. Limit theorems in one dimension\endheading

In this section, we will prove Theorem 3. We need to consider the first and second
moments of $W_t$. By duality,
$$\aligned EW_t=&E^{\eta}\sum_{x>0}\eta_t(x)=\sum_{x>0}P^{\eta}(\eta_t(x)=1)\\
=&\sum_{x>0}P^x(X_t\leq 0)=\sum_{y\leq 0<x}p_t(x,y)\\=&\sum_{n=1}^{\infty}
np_t(0,n)=E^0X_t^+.\endaligned\tag 4.1$$
Similarly,
$$\sum_{x>0}\big[P^{\eta}(\eta_t(x)=1)\big]^2=\sum_{x>0; u,v\leq 0}p_t(x,u)p_t(x,v)
=E^{(0,0)}\min(X_t^+,Y_t^+),\tag 4.2$$
where $Y_t$ is an independent copy of $X_t$.

For the covariances, we proceed as in Section 3, but this time with $f(x,y)=1_{\{x,y\leq 0\}}$ and
$g(x,y)=1_{\{x,y>0\}}$. For $x\neq y$, $(U-V)f(x,y)=p(x,y)[f(x,x)+f(y,y)-2f(x,y)]$, so for $x\neq y$,
$$\align (U-V)U(s)f(x,y)=& p(x,y)[U(s)f(x,x)+U(s)f(y,y)-2U(s)f(x,y)]\\=&p(x,y)[P^x(X_s\leq 0)-P^y(X_s
\leq 0)]^2\\=&p(x,y)[P^0(X_s\geq x)-P^0(X_s\geq y)]^2=\Delta_s(x,y),\endalign$$
where we define for all $x,y$,
$$\Delta_s(x,y)=\cases p(x,y)[P^0(x\leq X_s<y)]^2\quad&\text{if }x<y,\\
p(x,y)[P^0(y\leq X_s<x)]^2\quad&\text{if }y<x,\\0&\text{if }y=x.
\endcases$$
Then letting 
$$K(t)=-\sum_{x,y>0;x\neq y}Cov(\eta_t(x),\eta_t(y)),$$
we see that
$$\aligned K(t)=&\sum_{x,y>0;x\neq y}[U(t)-V(t)]f(x,y)
=\sum_{x,y>0;x\neq y}\int_0^tV(t-s)(U-V)U(s)f(x,y)ds\\=&\int_0^t\sum_{x\neq y}g(x,y)V(t-s)
\Delta_s(x,y)ds=\int_0^t\sum_{x,y}\Delta_s(x,y)V(t-s)g(x,y)ds\\\leq
&\int_0^t\sum_{x,y}\Delta_s(x,y)U(t-s)g(x,y)ds\\=&\int_0^t
\sum_{x,y}\Delta_s(x,y)P^0(X_{t-s}<x)P^0(X_{t-s}<y)ds.\endaligned\tag 4.3$$

Let $\rho(t,x)=P^0(X_t<x)$. Since $\Delta_s(x,y)=\Delta_s(y,x)=\Delta_s(1-x,1-y)$, we may continue (4.3)
by writing
$$\aligned K(t)\leq&\frac
12\int_0^t\sum_{x,y}\Delta_s(x,y)[\rho(t-s,x)\rho(t-s,y)+\rho(t-s,1-x)\rho(t-s,1-y)]ds
\\=&\int_0^t\sum_{x<y}p(x,y)[P^0(x\leq X_s<y)]^2\gamma(t-s,x,y)ds,\endaligned\tag 4.4$$
where
$$\gamma(t,x,y)=\rho(t,x)\rho(t,y)+\rho(t,1-x)\rho(t,1-y).$$
Note that since $\rho(t,x)+\rho(t,1-x)= 1$, it follows that $\gamma(t,x,y)\leq 1$. Now let
$$\Gamma(t,n,u,v)=\sum_{x:x\leq u,v<x+n}\gamma(t,x,x+n)\leq (n-|v-u|)^+.\tag 4.5$$
Then using the symmetry and translation invariance of $p_s(x,y)$,
$$\aligned K(t)\leq& \int_0^t\sum_{n=1}^{\infty}p(n)\sum_{u,v}p_s(0,u)p_s(0,v)\Gamma(t-s,n,u,v)ds
\\=&\int_0^t\sum_{n=1}^{\infty}p(n)E^0\Gamma(t-s,n,X_s,X_s-X_{2s})ds\\
\leq&\int_0^{t}\sum_{n=1}^{\infty}p(n)E^0(n-|X_{2s}|)^+ds.\endaligned\tag 4.6$$

Now assume that $\sigma^2=\sum_nn^2p(n)<\infty$, and write
$$E^0(n-|X_{2s}|)^+=\sum_{k=-n}^n(n-|k|)p_{2s}(0,k)\leq n^2p_{2s}(0,0),$$
where the inequality comes from
$$[p_{2s}(0,k)]^2=\bigg[\sum_jp_s(0,j)p_s(j,k)\bigg]^2\leq\sum_j[p_s(0,j)]^2\sum_j[p_s(j,k)]^2
=[p_{2s}(0,0)]^2.$$
Therefore,
$$\limsup_{t\rightarrow\infty}\frac 1{\sqrt t}\int_0^{t}\sum_{n>N}p(n)E^0(n-|X_{2s}|)^+ds\leq \frac
1{\sigma\sqrt{\pi}}\sum_{n>N}n^2p(n).\tag 4.7$$ 
by the local central limit theorem. 
For the terms corresponding to small $n$, we have given up too much in using the inequality
in (4.5). To handle these terms, first note that for fixed $k$, conditionally on $|X_{2s}|=k$,
$X_{s}/\sqrt s$ is asymptotically normal with mean 0 and variance $\sigma^2/2$.
Now take fixed $k,n$ with $0\leq k\leq n$.  Then

$$\align E^0[\Gamma(t-s,n,X_s,X_s-X_{2s})\mid X_{2s}=k]=& E^0\bigg[\sum_{X_s-n<x\leq
X_s-k}\gamma(t-s,x,x+n)\bigg| X_{2s}=k\bigg]\\=&E^0\bigg[\sum_{y=k}^{n-1}\gamma(t-s,X_s-y,X_s-y+n)
\bigg| X_{2s}=k\bigg].\endalign$$
Write
$$\align E^0[\gamma(t-s,X_s-y,&X_s-y+n)\mid X_{2s}=k]=P^0(Y_{t-s}<X_s-y,Z_{t-s}<X_s-y+n\mid X_{2s}=k)
\\&+P^0(Y_{t-s}<1-X_s+y,Z_{t-s}<1-X_s+y-n\mid X_{2s}=k),\endalign$$
where $X_t,Y_t,Z_t$ are independendent copies of the random walk. If $s,t\rightarrow\infty$ with
$s/r\rightarrow r\in (0,1)$, then the above expression converges to 
$$\gather P\bigg(N_2\leq\sqrt{\frac r{2(1-r)}}N_1,N_3\leq\sqrt{\frac r{2(1-r)}}N_1\bigg)\\+
P\bigg(N_2\leq-\sqrt{\frac r{2(1-r)}}N_1,N_3\leq-\sqrt{\frac r{2(1-r)}}N_1\bigg),\endgather$$
where $N_1,N_2,N_3$ are independent normally distributed random variables with mean
zero and variance 1. Call this expression $h(r)$.  Note that $h(r)<1$ for $r<1$,
so that
$$H=\int_0^1\frac {h(r)}{2\sqrt r}dr<1.$$

Passing to the limit in (4.6), we see that
$$\limsup_{t\rightarrow\infty}\frac {K(t)}{\sqrt t}\leq \frac
H{\sigma\sqrt{\pi}}\sum_{n=1}^{N}n^2p(n)+\frac
1{\sigma\sqrt{\pi}}\sum_{n>N}n^2p(n).$$
Taking $N\rightarrow\infty$ gives
$$\limsup_{t\rightarrow\infty}\frac {K(t)}{\sqrt t}\leq \frac{H\sigma}{2\sqrt{\pi}}.$$
Recalling from (4.1) and
(4.2) that
$$Var(W_t)=E^0X_t^+-E^0\min(X_t^+,Y_t^+)-K(t),$$
we see that
$$\limsup_{t\rightarrow\infty}t^{-1/2}Var(W_t)\leq \frac{\sigma}{2\sqrt{\pi}}$$
and
$$\liminf_{t\rightarrow\infty}t^{-1/2}Var(W_t)\geq\frac {(1-H)\sigma}{2\sqrt{\pi}}>0$$
So, Theorem 3 follows from Proposition 2.2(b) and the strong Rayleigh property of
$\{\eta_t(x), x>0\}$.

\heading 5. Stability and the Symmetric Exclusion Process\endheading

In this section, we present an elementary proof of the basic fact needed to prove preservation
of stability by the symmetric exclusion process. A similar proof was obtained independently
by Borcea and Br\"and\'en (2008),

A  polynomial $Q(z_1,...,z_n)$
with complex coefficients is said to be stable if $Q(z_1,...,z_n)\neq 0$ whenever $\Im(z_i)>0$ for each $i$.
The key result needed to show that the generating polynomial of the distribution a symmetric exclusion process
at time $t$ is stable whenever this is the case at time 0 is the following.

\proclaim {Theorem 5.1} Suppose the multi-affine polynomial $Q$ is stable. Then so is the
polynomial $Q_p$ for $0\leq p\leq 1$, where 
$$Q_p(z_1,...,z_n)=pQ(z_1,...,z_n)+(1-p)Q(z_2,z_1,z_3,...,z_n).$$
\endproclaim

The proof of Theorem 5.1 is based on the following characterization of stability for multi-affine
polynomials in two variables.

\proclaim{Proposition 5.2} Suppose $h(z,w)=a+bz+cw+dzw$, where $a,b,c,d$ are complex, and not all zero.
Then $h$ is stable if and only if
$$\Re(b\overline c-a\overline d)\geq |bc-ad|,\ \Im(a\overline b)\geq 0,\ \Im(a\overline c)\geq 0,
\ \Im(b\overline d)\geq 0,\ \Im(c\overline d)\geq 0.\tag 5.1$$
\endproclaim
\demo {Proof} If $b=c=d=0$, the $h$ is automatically stable, since then $a\neq 0$. If $d=0,b\neq 0$,
then $h(z,w)=0$ iff 
$$z=-\frac {(a+cw)\overline b}{|b|^2},$$
so stability is equivalent to $\Im(w)>0\Rightarrow \Im(z)\leq 0$, where
$$\Im(z)=-\frac{\Im(a\overline b)+\Re(c\overline b)\Im (w)+\Im(c\overline
b)\Re(w)}{|b|^2}.$$
Therefore, since $\Re(w)$ is arbitrary, stability is equivalent to $c\overline b\geq 0$ and $\Im(a\overline
b)\geq 0$, and these imply $\Im (a\overline c)=b\overline c\Im(a\overline b)/|b|^2\geq 0$ in this case. 

So, we may now assume that $d\neq 0$. Solving for $z$, we see that $h(z,w)=0$ iff 
$$b+dw=0,\quad a+cw=0\tag 5.2$$
or
$$b+dw\neq 0,\quad z=-\frac{a+cw}{b+dw}=-\frac{a\overline b+c\overline bw+a\overline d\overline w
+c\overline d|w|^2}{|b+dw|^2}.\tag 5.3$$
Case (5.2) occurs iff
$$w=-\frac bd=-\frac{b\overline d}{|d|^2},\quad bc=ad.\tag 5.4$$
Therefore, if $bc\neq ad$, stability is equivalent to the statement 
$$\Im (w)\geq0\Rightarrow
\Im(a\overline b)+\Re(w)\Im(c\overline b+a\overline d)+\Im(w)\Re(c\overline b-a\overline
d)+|w|^2\Im(c\overline d)\geq0,\tag 5.5$$
while if $bc=ad$, stability is equivalent to this, together with $\Im(b\overline d)\geq 0$.

Letting $w$ be real, we see that for (5.5) to hold, we need $$\Im(a\overline b)\geq 0,\quad
\Im(c\overline d)\geq 0, \quad
[\Im(c\overline b+a\overline d)]^2\leq 4\Im(a\overline b)\Im(c\overline d).$$
Note that
$$4\Im(a\overline b)\Im(c\overline d)-[\Im(c\overline b+a\overline d)]^2=
[\Re(c\overline b-a\overline d)]^2-|bc-ad|^2.$$
Minimizing the expression in (5.5) over $\Re(w)$, we see that we also need $\Re(c\overline b-a\overline d)\geq
0$ if $\Im (c\overline d)=0$, and 
$$\Im(a\overline b)-\frac{[\Im(c\overline b+a\overline d)]^2}{4\Im(c\overline d)}+\Re(c\overline b-a
\overline d)t+\Im (c\overline d)t^2\geq 0,\quad t\geq 0\tag 5.6$$
if $\Im (c\overline d)>0$. Since the discriminant of this quadratic is
$$[\Re(c\overline b-a\overline d)]^2+[\Im(c\overline b+a\overline d)]^2
-4\Im(a\overline b)\Im(c\overline  d)=|bc-ad|^2\geq 0,$$
if $\Im(c\overline d)>0$, (5.6) is always true if $bc=ad$, while if $bc\neq ad$, (5.6) is equivalent to 
$\Re(c\overline b-a\overline d)\geq 0$. Putting these observations together, noting that
stability is not changed if the roles of $b$ and $c$ reversed, completes the proof.
\enddemo

\demo{Proof of Theorem 5.1} We need to show that $Q_p(z_1,...,z_n)\neq 0$ whenever $\Im(z_i)>0$ for
$i=1,...,n$. To to so, fix $z_3,...z_n$ with $\Im(z_i)>0$ for $i=3,...,n$, and write
$$h(z,w)=Q(z,w,z_3,...,z_n).$$
Then $h$ is of the form considered in Proposition 5.2, and we must show that if (5.1) holds
for a given $a,b,c,d$, then it also holds with $b$ and $c$ replaced by
$$b(p)=pb+(1-p)c,\quad c(p)=pc+(1-p)b.$$
This follows from
$$\gather\Im(a\overline{b(p)})=p\Im(a\overline b)+(1-p)\Im(a\overline c),\quad
\Im(a\overline{c(p)})=p\Im(a\overline c)+(1-p)\Im(a\overline b),\\
\Im(b(p)\overline{d})=p\Im(b\overline d)+(1-p)\Im(c\overline d),\quad
\Im(c(p)\overline{d})=p\Im(c\overline d)+(1-p)\Im(b\overline d),
\endgather$$
and
$$\gather\Re(b(p)\overline{c(p)}-a\overline d)=\Re(b\overline{c}-a\overline d)+p(1-p)|b-c|^2,\\
b(p)c(p)-ad=(bc-ad)+p(1-p)(b-c)^2,\endgather$$
so that
$$\Re(b(p)\overline{c(p)}-a\overline d)-|b(p)c(p)-ad|\geq \Re(b\overline{c}-a\overline d)-|bc-ad|$$
by the triangle inequality.\enddemo

\proclaim{Corollary 5.3} Suppose the generating polynomial of the initial distribution of a symmetric
exclusion process on a finite set $S$ is stable. Then the same is true at time $t>0$. \endproclaim
\demo{Proof} View the process in terms of stirrings. In other words, for each pair $x,y\in S$
interchange $\eta(x)$ and $\eta(y)$ at Poisson times at a certain rate. If stirrings are applied
at only one pair of sites, the generating polynomial of the distribution at time $t$ is of the
form $Q_p$ given in the statement of Theorem 5, where $Q$ is the generating polynomial of the
initial distribution. For a general exclusion process, the distribution at time $t$ can be
obtained as a limit of that obtained by successively applying stirrings at different pairs
of sites.
\enddemo

\bigskip
\centerline{\bf REFERENCES}
\bigskip
\ref\by R. Arratia\paper The motion of a tagged particle in the simple symmetric exclusion
process on $Z$\jour Ann. Probab.\vol 11\yr 1983\pages 362--373\endref

\ref\by A. D. Barbour, L. Holst and S. Janson\book Poisson Approximation\publ Oxford University
Press\yr 1992\endref

\ref\by J. Borcea and P. Br\"and\'en\paper Linear operators preserving stability\yr 2008\endref

\ref\by J. Borcea, P. Br\"and\'en and T. M. Liggett\paper Negative dependence and the
geometry of polynomials\yr 2008\endref

\ref\by P. Br\"and\'en \paper Polynomials with the half-plane property and matroid
theory\jour Adv. Math.\yr 2007\vol 216\pages 302--320\endref

\ref\by D. I. Cartwright and W. Woess\paper Infinite graphs with nonconstant Dirichlet
finite harmonic functions\jour SIAM J. Disc. Math.\vol 5\yr 1992\pages 380--385\endref

\ref\by A. De Masi and P. A. Ferrari\paper Flux fluctuations in the one dimensional nearest
neighbors symmetric simple exclusion process\jour J. Stat. Phys.\vol 107\yr 2002\pages
677-683\endref

\ref \by M. D. Jara and C. Landim\paper Nonequilibrium central limit theorem for a
tagged particle in symmetric simple exclusion\jour Ann. I. H. Poincar\'e\vol 42\pages
567--577\yr 2006\endref

\ref\by C. Kipnis and S. R. S. Varadhan\paper Central limit theorem for additive functionals
of reversible Markov processes and applications to simple exclusions\jour Comm. Math. Phys.
\vol 104\yr 1986\pages 1--19\endref

\ref \by T. M. Liggett\book Interacting Particle Systems\publ Springer-Verlag
\publaddr New York\yr 1985\endref

\ref\by C. M. Newman\paper Asymptotic independence and limit theorems for positively and
negatively dependent random variables\inbook Inequalities in Statistics and Probability
\publ IMS\yr 1984\pages 127--140\endref

\ref\by M. Peligrad and S. Sethuraman\paper On fractional Brownian motion limits in one
dimensional nearest-neighbor symmetric simple exclusion\yr 2008\endref

\ref\by R. Pemantle \paper Towards a theory of negative
dependence\jour J. Math. Phys.\vol 41
\pages 1371--1390\yr 2000\endref

\ref\by G. G. Roussas\paper Asymptotic normality of random fields of positively or negatively
associated processes\jour J. Mult. Anal.\vol 50\yr 1994\pages 152--173\endref

\bigskip

\noindent Department of Mathematics

\noindent University of California, Los Angeles

\noindent 405 Hilgard Ave.

\noindent Los Angeles CA 90095
\bigskip
\noindent email: tml\@math.ucla.edu

\noindent URL:  http://www.math.ucla.edu/\~\ tml/
\end